\input amstex
\magnification 1200
\documentstyle{amsppt}

\pagewidth{16.2 truecm}
\pageheight{22.5 truecm}
\hcorrection{-0.3 truecm}
\vcorrection{0.6 truecm}

\parskip 5pt

\NoBlackBoxes
\nologo

\define\a{\alpha}

\define\diz{(1-|z|^2)}
\define\eit{e^{i\theta}}

\define\disank{(1-|a_{n,k}|^2)}

\define\sak{\{a_k\}}

\redefine\phi{\varphi}

\topmatter
\title
Sampling sets for Hardy spaces of the disk
\endtitle

\author
Pascal J. Thomas
\endauthor

\affil
Universit\'e Paul Sabatier
\endaffil

\address
{Pascal J. Thomas:
Laboratoire Emile Picard
Universit\'e Paul Sabatier,
118 route de Narbonne,
31062 TOULOUSE CEDEX,
France; pthomas\@cict.fr
}
\endaddress

\abstract{We propose two possible definitions for the
notion of a sampling sequence (or set) for Hardy spaces of
the disk. The first one is inspired by recent work of Bruna,
Nicolau, and \O yma about interpolating sequences in the same
spaces, and it yields sampling sets which do not depend on
the value of $p$ and correspond to the result proved
for bounded functions ($p=\infty$) by Brown, Shields and Zeller.
The second notion, while formally closer to the one used
for weighted Bergman spaces, is shown to lead
to trivial situations only, but raises a possibly
interesting problem.}
\endabstract

%\subjclass{}
%\endsubjclass

\endtopmatter
\document
\subheading{\S 1. Sampling via the Bruna-Nicolau-\O yma function}

Let $\Bbb D$ be the unit disk in the complex plane.
Recall that for any $0<p\leq \infty$
the Hardy space $H^p ( \Bbb D)$ is the
set of holomorphic functions $f$ such that
$$
\| f \|_p := \left( \sup_{r<1}
\int_0^{2\pi} |f(re^{i\theta})|^p \, \frac{d \theta}{2 \pi}
\right)^{\frac1p} < \infty ,
$$
where the integral is replaced by a supremum in the case
$p=\infty$, and that for $p \geq 1$, $\| \cdot \|_p$ is
a norm.

In accordance with many previous works (e.g. \cite{La},
\cite{Se1}, \cite{Se2}) we would
like to say that a subset $a$ of the unit disk is of {\it sampling \/}
for the space $H^p ( \Bbb D)$ when the values of a function
$f \in H^p ( \Bbb D)$, restricted to this set, determine
the function uniquely, and when we can establish some
inequalities between the $H^p$-norm and an appropriate norm
on the space of functions on the subset $a$.
Usually, this is interesting only when the subset $a$
is a discrete sequence of points, however, that hypothesis
will not be necessary for the first part of this paper.
We proceed to give a more specific definition.

The {\it Stolz angle\/} with vertex at $e^{i\theta}$ is
 $\Gamma_\alpha (e^{i\theta}) :=
\{ z \in \Bbb D : \frac{|e^{i\theta}-z|}{1-|z|} < 1 + \alpha \}$.

The {\it nontangential maximal function\/} is
$$
Mf (e^{i\theta}) := \sup_{z \in \Gamma_\alpha (e^{i\theta}) } |f(z)|.
$$
For any $0<p\leq \infty$, for any choice of
$\alpha >0$, we have $\| Mf \|_p \leq C_{p, \alpha} \| f \|_p$,
where, for functions defined on the unit circle,
$\| \cdot \|_p$ stands for the usual norm in the space
$L^p (\frac{d \theta}{2 \pi})$ (\cite{Du}, \cite{Ga, Theorem II.3.1,
p.~57}).

Now, following Bruna-Nicolau-\O yma \cite{Br-Ni-\O y}, let
$$
M_a (f) (e^{i\theta}) :=
M_{a,\infty} (f) (e^{i\theta})
:= \sup_{\Gamma_\alpha (e^{i\theta}) \cap a} |f|
\leq Mf (e^{i\theta}).
$$

From the above it follows that $\| M_a (f) \|_p \leq C_{p, \alpha} \| f \|_p$.
We will call the set $a$ sampling if the two norms are
actually equivalent.

\proclaim{Definition 1}

We say that the set $a$ is sampling for $ H^p ( \Bbb D)$
iff there exists a constant $C>0$ such that for any
$f \in H^p ( \Bbb D)$, $\| M_a (f) \|_p \geq C \| f \|_p$.
\endproclaim

In the case where $p=\infty$, this simply says that
$\sup_a |f| \geq C \sup_{\Bbb D} |f|$, and by taking powers
of $f$ we see that $\sup_a |f| = \sup_{\Bbb D} |f|$. This case of
the problem was solved by Brown, Shields and Zeller
\cite{Br-Sh-Ze, Th.~3, (iii)-(iv)}.
The main positive result of this note is that, with the appropriate
definition above, this can be extended to all $p>0$.

\proclaim{Definition}
We say that a point of the circle $e^{i\theta}$ is a {\it nontangential
limit point\/} of the set $a$ iff it is in the closure of
$\Gamma_\alpha (e^{i\theta}) \cap a$ for some $\alpha > 0$. In this
case, we write $e^{i\theta} \in NT(a)$.
\endproclaim

Denote by $\lambda_1$ the $1$-dimensional Lebesgue measure
on the unit circle.

\proclaim{Theorem 1}

$a$ is sampling for $ H^p ( \Bbb D)$ if and only if
$\lambda_1$-almost every point  $e^{i\theta} \in \partial \Bbb D$
 is a nontangential limit point of the set $a$.
\endproclaim

\demo{Proof}
This proof is essentially the same as in \cite{Br-Sh-Ze}, so
we shall keep it brief. We say that
a function $f$ defined on the disk admits a {\it nontangential
limit\/} at the point $e^{i\theta}$ iff
$\lim_{z \to e^{i\theta}, z \in \Gamma_\alpha (e^{i\theta}) } f(z) =:
f^* (e^{i\theta}) $ exists and is finite.
For for any choice of
$\alpha$, any $f \in H^p ( \Bbb D)$ admits a nontangential
limit at almost every $e^{i\theta}$
and $\| f^* \|_p =\| f \|_p $ \cite{Du}, \cite{Ga}.

Thus for almost every $e^{i\theta} \in NT(a)$,
$$
M_a (f) \geq
\lim_{z \to e^{i\theta}, z \in \Gamma_\alpha (e^{i\theta}) \cap a} | f(z) |
= | f^* (e^{i\theta})| ,
$$
so that if $\lambda_1 (\partial \Bbb D \setminus NT(a))=0$,
$\| M_a (f) \|_p \geq \| f^* \|_p$, q.e.d.

Conversely, if $\lambda_1 (\partial \Bbb D \setminus NT(a) )>0$,
there is an
integer $N$ and $A$, a compact set of positive measure, such that
for all $e^{i\theta} \in A$,
$\Gamma_\alpha (e^{i\theta}) \cap a \cap \{ |z| \geq 1- 1/N \} =
\emptyset$. Consider the outer function
$$
\omega_A (z) := \exp \left\{ - \int_0^{2\pi}
\frac{z+e^{i\theta}}{z-\eit} (1 - \Bbb 1_A (\eit))
\frac{d \theta}{2 \pi} \right\},
$$
where $\Bbb 1_A$ is the indicator function of $A$. Then
$-\log | \omega_A |$ is the harmonic measure of the set
$\partial \Bbb D \setminus A$,
and classical estimates (e.g.\ \cite{Ga, ex.~3, p.~41})
show that for all
$z \notin \cup_{\eit \in A} \Gamma_\a (\eit)$,
$-\log | \omega_A | \geq c_\a > 0$.

Consider the sequence
of functions $f_n(z) := z^n \omega_A (z)^n$ : we have
$\| f_n \|_p = \| f_n^* \|_p $ $\geq \frac{\lambda_1 (A)}{2 \pi} > 0$,
while $M_a(f_n)(\eit) \leq 1$ for all $n$,
and $\lim_{n \to \infty} M_a(f_n)(\eit) = 0$ for all
$\eit \in A$ because $(1-1/N)^n \to 0$, and for all
$\eit \notin A$ because $\exp(-n c_\a) \to 0$. So
$\| M_a (f_n) \|_p \to 0$, and the sampling inequality cannot hold.

\enddemo

\subheading{\S 2. Attempt at a classical definition}

The Bruna-Nicolau-\O yma function was introduced to deal with
problems of interpolation of the type studied in \cite{Sh-Sh}.
This involved a norm in $L^p(\mu)$, where for any $E \subset a$,
$\mu (E) := \sum_{z\in E} \diz $.

We can then state a similar sampling problem:

\proclaim{Definition 2}

We say that the set $a$ $ H^p ( \Bbb D)${\it-thick}
iff there exists a constant $C>0$ such that for any
$f \in H^p ( \Bbb D)$,
$$
\| f \|_{L^p(\mu)} =
\left( \sum_{z \in a} \diz |f(z)|^p \right)^{\frac1p} \geq C \| f \|_p .
$$
\endproclaim

This says that the measure $\mu$ is {\it dominating\/}
in the sense of \cite{Lu}.

Note that the sum on the left hand side may be infinite;
 in fact it always will be, which makes
 for a rather uninteresting notion of sampling.

\proclaim{Theorem 2}
\roster
\item
If $\lambda_1 (NT(a)) > 0$, then
$H^p ( \Bbb D) \cap L^p(\mu) = \{ 0 \}$
(i.e. if $f\in H^p ( \Bbb D)$ and
$\sum_{z \in a} \diz |f(z)|^p < \infty$, then $f=0$).
\item
The set $a$ is $H^p$-thick if and only if
$H^p ( \Bbb D) \cap L^p(\mu) = \{ 0 \}$ .
\endroster
\endproclaim

In the case of the Bergman spaces studied in \cite{Se1},
\cite{Se2}, the definition of sampling which is given is
akin to Definition 2 (but, unlike it, is not
vacuous!). One can also construct an
analogue to the first definition, using a supremum on invariant balls
of fixed radius rather than Stolz angles.
It can be seen, using boundedness of the restriction
map which is part of the definition,
that the sequences under consideration in \cite{Se1, \S 7}
can only have a bounded number of points in each such ball, which
implies that the two notions are in this case equivalent.
The proof is similar to that of the following Lemma.

\proclaim{Lemma 1}

If $a$ is sampling for $ H^p ( \Bbb D)$, then $a$
is $ H^p ( \Bbb D)$-thick.
\endproclaim

\demo{Proof of Lemma 1}
Let
$$
M_{a,p} (f) (e^{i\theta}) :=
\left( \sum_{z \in \Gamma_\alpha (e^{i\theta}) \cap a} |f(z)|^p
\right)^{\frac1p} .
$$
Then
$M_{a,p} (f) (e^{i\theta}) \geq M_{a,\infty} (f) (e^{i\theta})$
for each $\theta$. But
$$
\int_0^{2\pi} M_{a,p} (f) (e^{i\theta})^p d \theta
=
\sum_{z \in \Gamma_\alpha (e^{i\theta}) \cap a} |f(z)|^p
\int_{\theta : z \in \Gamma_\alpha (e^{i\theta})} d \theta ,
$$
and the set $I_z := \{ \theta : z \in \Gamma_\alpha (e^{i\theta}) \} $
is an arc centered at $z/|z|$, with length a multiple of
$1 - |z|$, so the last sum is commensurate to the one in
Definition 2.
\enddemo

\demo{Proof of Theorem 2}
The following proof of \therosteritem1 is due to Bo Berndtsson.

Assume $f \in L^p(\mu)$. Then, by the proof of Lemma 1,
$M_{a,p} (f) $ must be finite almost everywhere on $\partial \Bbb D$.
At the points  $\eit \in NT(a)$, since the infinite sum converges,
$$
\lim_{z \in \Gamma_\alpha (e^{i\theta}) \cap a, z \to \eit} |f(z)|^p =0 .
$$
Suppose now that in addition $f \in H^p ( \Bbb D)$.
Since $f$ admits a nontangential limit almost everywhere,
that limit will be zero almost everywhere on $NT(a)$, a set of
positive measure, thus $f=0$ (see e.g.\
\cite{Ga, Th.~4.1, p.~64}).

As remarked after Definition 2, if
 $L^p(\mu) \cap H^p (\Bbb D)= \{ 0 \}$, then $a$ is
$H^p$-thick in a trivial way. To prove the converse implication
suppose $0 \neq f \in L^p(\mu) \cap H^p (\Bbb D)$. We shall
prove that $NT(a)$ is of full measure. This leads to a
contradiction in view of \therosteritem1.

Assume $NT(a)$ is not of full measure. We can then construct
the same sequence $f_n$ as in the proof of Theorem 1, and we have
$\lim_{n \to \infty} f_n (z) = 0$ for all $z \in a$. Let $f$ be a
function as in property ($F_p$); then the sequence $\{f f_n\}$ is
dominated by $|f|$, and we can apply Lebesgue's dominated convergence
theorem in $L^p (\mu)$, where for any $E \subset a$,
$\mu (E) := \sum_{z\in E} \diz $. So
$\lim_{n\to \infty} \sum_{z\in E} \diz |f f_n|^p =0$, while
$\| f f_n \|_p \geq \int_A |f^*| > 0$, since
$f^*$ cannot vanish on a set of positive measure.

\enddemo

\subheading{\S 3. A question, and an example}

\proclaim{Question}
Can we characterize explicitly the sets $a$
such that $L^p(\mu) \cap H^p (\Bbb D) \neq \{ 0 \}$ ?
\endproclaim

The sets under consideration will have to be
discrete sequences.
I see the property as a weaker analogue of the Blaschke property
(a sequence satisfying
the Blaschke property, being a zero-set for an $H^\infty$
function, automatically satisfies our condition).

The condition cannot be
about the mere growth of the number of points in the sequence as it
approaches the boundary, as is demonstrated by the following
example: take $\{b_n\}$ a sequence of points in the disk
satisfying the Blaschke condition, and
$\{ b_{n,k} , 1\leq k \leq q_n\}$, $q_n$ distinct points in
the disk $D(b_n , q_n^{-1} (1-|b_n|^2)^2)$. It is easy to
check that the Blaschke product with simple zeroes at each
of the $b_n$ is in $L^p(\mu) \cap H^p (\Bbb D)$
for any $p\geq 1$, but of course
$q_n$ can grow as fast as we please.

The condition does not
depend on $p$,
in keeping with the fact that zero-sets or interpolating
sequences for the Hardy spaces $H^p (\Bbb D)$ do not depend on $p$.

\proclaim{Lemma 2}

There exists
$H^p (\Bbb D) \cap L^p(\mu) \neq \{0\} $
if and only if
$H^\infty (\Bbb D) \cap L^1(\mu) \neq \{0\}$.
\endproclaim

\demo{Proof of Lemma 2}
This uses the same ideas as the canonical
factorization of functions in the Nevanlinna class.
Let $f \in H^p (\Bbb D) \cap L^p(\mu) \setminus \{0\} $.
Let $m(z)$ be a harmonic majorant of $|f|^p$ \cite{Ga, p.~50}.
Since $m(z) \geq 0$, if we
denote by $\tilde m$ its harmonic conjugate,
$1+m(z)+i\tilde m(z) = e^{H(z)}$, with
$e^{\frac1p \Re H(z)} \geq \max(1,|f(z)|)$. Then
$$
f_1 (z) := e^{- \frac1p  H(z)} f(z) \in H^\infty (\Bbb D) ,
\quad \|f_1 \|_\infty \leq 1 ,
$$
and $|f_1 (z)| \leq |f (z)|$, so $f_1 \in L^p(\mu)$.
If $p\leq 1$, $|f_1 (z)|^p \geq |f_1 (z)|$ and we may set
$g = f_1$.
Otherwise, we write $f_1 (z) = B (z) e^{h(z)}$,
where $B$ is a Blaschke product, and
set $g(z) =  B (z)^{[p]+1} e^{p h(z)}$; then $|g(z)| \leq |f_1 (z)|^p$,
and we're done.

The converse implication is easier.
\enddemo

\proclaim{Proposition 3}
There exists a discrete sequence, accumulating at one point
of the circle only, such that
$H^p (\Bbb D) \cap L^p(\mu) = \{0\} $.
\endproclaim

\demo{Proof}
Let $\{ p_n \}$ be chosen so that $2^{-n} p_n \to \infty$, and
take $\gamma_n := 2^n/2p_n$.
\newline
Set $\ell_n := p_n^{-1} |\log \gamma_n |^{-1/2}$,
and
$$
a_{n,k} := (1-2^{-n})^{\frac12} \exp {i k \ell_n} , \quad
0 \leq k \leq p_n .
$$
Let $g$ be a function as in Lemma 2;
to prove that $g=0$, it is enough to
show that
$\lim_{r \to 1^- } \int_0^{2\pi} \log |g(r\eit)| d \theta
= -\infty$
\cite{Ga, proof of Th.~4.1, p.~65}.  We will prove that
$|g|$ is small enough often enough on circles of radii
tending to $1$.

Set
$$
s_n := \sum_{k\leq p_n} \disank |g(a_{n,k})|
= 2^{-n} \sum_{k\leq p_n}  |g(a_{n,k})| .
$$
Since $\sum_n s_n < \infty$, $\lim_{n\to\infty} s_n =0$.
Let
$$
A_n = \{ k : 1 \leq k \leq p_n , |g(a_{n,k})| \geq \gamma_n \} .
$$
By Chebyshev's inequality, $\# A_n \leq 2^n s_n/\gamma_n \leq p_n /2$,
thus if we set
$A'_n := \{1, \dots , p_n\} \setminus A_n$,
$\#A'_n \geq p_n /2$.

Set
$ J(a_{n,k}) := ((k-\frac12) \ell_n , (k+\frac12) \ell_n)$. Then
$$ \lambda_1 \left( \cup_{k\in A'_n} J(a_{n,k}) \right)
\geq C p_n \ell_n .
$$
Since $g \in H^\infty (\Bbb D)$, $|g'(z)| \leq C \diz^{-1}$,
so that for $\theta \in J(a_{n,k})$,
$$
\left| g( (1-2^{-n})^{\frac12} \eit ) - g (a_{n,k}) \right|
\leq C {2^{n}}{\ell_n} ,
$$
thus when $k \in A'_n$,
$$
 \log | g( (1-2^{-n})^{\frac12} \eit ) | \leq
\log (\gamma_n + C 2^n \ell_n )  \leq C + \log \gamma_n .
$$

The integral is now estimated by
$p_n \ell_n \log \gamma_n \leq - C |\log \gamma_n |^{1/2} \to -\infty$,
and the arc subtended by the points $a_{n,k}$ for a fixed $n$
is of length $p_n \ell_n = |\log \gamma_n |^{-1/2} \to 0$ as
$n \to \infty$, so $\sak$ only accumulates at $1$.

\enddemo

{\bf Acknowledgements.}
This paper arose as a by-product of my collaboration with
Xavier Massaneda, with whom it was much discussed.
The key observation for the proof of Theorem 2 is
due to Bo Berndtsson, who also raised the possibility
that sampling in the Hardy spaces could be an
interesting question.
Also I should thank Alex Heinis for bringing
\cite{Br-Sh-Ze} to my attention (while he was studying
under Jan Wiegerinck's direction).

\refstyle{A}

\enddocument